\PassOptionsToPackage{unicode}{hyperref}
\PassOptionsToPackage{hyphens}{url}
\PassOptionsToPackage{dvipsnames,svgnames,x11names}{xcolor}
\documentclass[
  american,
  11pt]{article}
\usepackage{xcolor}
\usepackage[margin=1in]{geometry}
\usepackage{amsmath,amssymb}
\usepackage{iftex}
\ifPDFTeX
  \usepackage[T1]{fontenc}
  \usepackage[utf8]{inputenc}
  \usepackage{textcomp} 
\else 
  \usepackage{fontspec}
  \defaultfontfeatures{Scale=MatchLowercase}
  \defaultfontfeatures[\rmfamily]{Ligatures=TeX,Scale=1}
\fi
\usepackage{lmodern}
\ifPDFTeX\else
\fi
\IfFileExists{upquote.sty}{\usepackage{upquote}}{}
\IfFileExists{microtype.sty}{
  \usepackage[]{microtype}
  \UseMicrotypeSet[protrusion]{basicmath} 
}{}
\makeatletter
\@ifundefined{KOMAClassName}{
  \IfFileExists{parskip.sty}{%
    \usepackage{parskip}
  }{
    \setlength{\parindent}{0pt}
    \setlength{\parskip}{6pt plus 2pt minus 1pt}}
}{
  \KOMAoptions{parskip=half}}
\makeatother
\usepackage{color}
\usepackage{fancyvrb}

\DefineVerbatimEnvironment{Highlighting}{Verbatim}{commandchars=\\\{\}}
\newenvironment{Shaded}{}{}

\newcommand{\NormalTok}[1]{#1}

\usepackage{longtable,booktabs,array}
\usepackage{calc} 
\usepackage{etoolbox}
\makeatletter
\patchcmd\longtable{\par}{\if@noskipsec\mbox{}\fi\par}{}{}
\makeatother
\IfFileExists{footnotehyper.sty}{\usepackage{footnotehyper}}{\usepackage{footnote}}
\makesavenoteenv{longtable}
\ifLuaTeX
\usepackage[bidi=basic,shorthands=off,]{babel}
\else
\usepackage[bidi=default,shorthands=off,]{babel}
\fi
\ifLuaTeX
  \usepackage{selnolig} 
\fi
\providecommand{\tightlist}{%
  \setlength{\itemsep}{0pt}\setlength{\parskip}{0pt}}
\usepackage{bookmark}
\IfFileExists{xurl.sty}{\usepackage{xurl}}{} 
\hypersetup{
  pdftitle={Uniformity without Projective Consistency: An Exact Counterexample for a Nested Binary Term Grammar},
  pdfauthor={Ivan Khalamendyk},
  pdfsubject={A constructive counterexample to projective consistency for uniform linear-extension measures},
  pdfkeywords={partially ordered set, linear extension, uniform measure, projective consistency, exact integer enumeration, central measure},
  pdflang={en-US},
  colorlinks=true,
  linkcolor={Maroon},
  filecolor={Maroon},
  citecolor={Blue},
  urlcolor={Blue},
  pdfcreator={LaTeX via pandoc}}

\title{Uniformity without Projective Consistency: An Exact
Counterexample for a Nested Binary Term Grammar}
\author{Ivan Khalamendyk\\[0.25em]\small Independent researcher}
\date{10 August 2026}

\begin{document}
\maketitle

\begin{abstract}

Consider the nested finite sets of ordered binary terms

\[
T_0=\{L\},
\qquad
T_{r+1}=\{L\}\cup\{N(a,b):a,b\in T_r\},
\]

and the partial orders on the events \(E_r=T_r\setminus\{L\}\) in which
an event \(N(a,b)\) may appear only after each of its nonleaf children.
At each level \(r\), let \(\mu_r\) be the uniform measure on the set of
linear extensions of this order. The natural map \(\rho_{43}\) deletes
all new events from a level-\(4\) linear extension while preserving the
relative order of the level-\(3\) events.

It is proved that the pushforward \((\rho_{43})_*\mu_4\) is not equal to
\(\mu_3\). Two explicit orders on the \(25\) events at level \(3\) are
constructed. For a prefix of length \(i\), the number of new events
already released has the closed form \((i+1)^2-b_i^2\), where \(b_i\) is
the number of terms from \(T_2\) in the prefix, including the leaf. The
release profile of the depth-priority order dominates that of the level
order pointwise and is strictly larger for \(3\le i\le15\). An explicit
injection between admissible interleavings converts this dominance into
a strict analytic proof that the two fibers have different
cardinalities.

Two algorithmically distinct exact counts reproduce both fiber
cardinalities. Each has \(1557\) decimal digits, and their reduced ratio
is

\[
\frac{614690215260160000}{479048686862260621}
\approx1.2831476885707443.
\]

The analogous restrictions through level \(3\) are consistent; thus
\(4\to3\) is the first failure in this grammar. The result is a
constructive counterexample specifically for the stated grammar, the
uniform family of measures, and the specified restriction map. It does
not exclude other weights, grammars, state spaces, or coarse-graining
rules.

\end{abstract}



\section{1. Problem statement}\label{problem-statement}

The uniform measure on a finite set is always unambiguously defined. A
family of uniform measures on nested finite spaces, however, need not be
consistent. If a larger space projects onto a smaller one, points in the
smaller space receive weights proportional to the numbers of their
extensions. Uniformity is preserved only when all fibers of the map have
the same cardinality. This general phenomenon is elementary: in the
three-element poset with the single relation \(a<x\), the three linear
extensions \(abx,axb,bax\) restrict after deletion of \(x\) to
\(ab,ab,ba\), producing probabilities \(2/3\) and \(1/3\) rather than a
uniform pair.

This elementary fact becomes nontrivial when the points are linear
extensions of rapidly growing partial orders. We study one completely
specified grammar. At every level, it produces a finite set of syntactic
events and a prerequisite order. The random object is an admissible
sequential listing of all events, with all such listings at a fixed
level assigned equal weight.

The main question is:

\begin{quote}
Does the distribution of old events obtained by deleting the new events
from a uniform level-\(4\) order coincide with the uniform distribution
on level-\(3\) orders?
\end{quote}

The answer is no. It can be proved without listing all level-\(3\)
orders and without computing the total number of level-\(4\) orders. It
is enough to find two lower-level orders with different numbers of
upper-level extensions. Below, the sign of this difference is proved
analytically, and its exact magnitude is computed in two different ways.

The general theory of random linear extensions, central measures,
order-invariant measures, and consistent random orderings is established
{[}1-8{]}. The contribution here is narrower: delayed failure in an
explicit recursive family that is consistent through the preceding
levels, a particular deletion map from level \(4\) to level \(3\), two
explicit fibers, a structural release-profile proof that their sizes
differ, and a reproducible exact count.

\section{2. Nested grammar and partial
orders}\label{nested-grammar-and-partial-orders}

\subsection{2.1. Terms}\label{terms}

Let \(L\) be a distinguished leaf and let \(N(\cdot,\cdot)\) be an
injective constructor for ordered pairs. Define recursively

\[
T_0=\{L\},
\qquad
T_{r+1}=\{L\}\cup\{N(a,b):a,b\in T_r\}.
\qquad\text{(1)}
\]

Injectivity means

\[
N(a,b)=N(c,d)
\quad\Longleftrightarrow\quad
a=c\ \text{and}\ b=d,
\qquad\text{(2)}
\]

and no term \(N(a,b)\) equals \(L\). The constructor is ordered: in
general, \(N(a,b)\) and \(N(b,a)\) are distinct. Repeated construction
of the same term does not create a new copy, because \(T_r\) is a set.

We have \(T_r\subseteq T_{r+1}\). Indeed, if \(N(a,b)\in T_r\), then
\(a,b\in T_{r-1}\subseteq T_r\), so the same term belongs to
\(T_{r+1}\). The numbers of terms satisfy the known recurrence

\[
t_0=1,
\qquad
t_{r+1}=1+t_r^2,
\qquad
t_r=|T_r|.
\qquad\text{(3)}
\]

Consequently,

\[
|T_0|,|T_1|,|T_2|,|T_3|,|T_4|
=1,2,5,26,677.
\qquad\text{(4)}
\]

The recurrence (3) and its initial values are known in the combinatorial
literature {[}9,10{]}; no novelty is claimed for them here.

\subsection{2.2. Events and
prerequisites}\label{events-and-prerequisites}

Define the event set at level \(r\) by

\[
E_r=T_r\setminus\{L\}.
\qquad\text{(5)}
\]

For an event \(t=N(a,b)\), its immediate prerequisites are

\[
\operatorname{Pred}(t)=\{a,b\}\setminus\{L\}.
\qquad\text{(6)}
\]

A repeated child gives a single prerequisite:
\(\operatorname{Pred}(N(a,a))=\{a\}\) for \(a\ne L\). Events are not
consumed, there are no conflicts between them, and no additional
dependency edges are introduced.

Let \(P_r=(E_r,\preceq_r)\) be the partial order obtained as the
reflexive-transitive closure of the relations

\[
x\prec_r N(a,b)
\quad\text{for every}\quad
x\in\operatorname{Pred}(N(a,b)).
\qquad\text{(7)}
\]

The height of a term,

\[
h(L)=0,
\qquad
h(N(a,b))=1+\max\{h(a),h(b)\},
\qquad\text{(8)}
\]

increases strictly along every prerequisite edge. Hence there are no
cycles, and \(P_r\) is indeed a finite partially ordered set.

\subsection{2.3. Linear extensions and uniform
measures}\label{linear-extensions-and-uniform-measures}

A linear extension of \(P_r\) is a sequence

\[
\sigma=(\sigma_1,\ldots,\sigma_{|E_r|})
\]

containing all events in \(E_r\), in which \(x\prec_r y\) implies that
\(x\) appears before \(y\). Denote the set of these sequences by
\(\operatorname{Lin}(P_r)\) and its cardinality by \(e(P_r)\). At each
level, define the uniform measure

\[
\mu_r(\{\sigma\})=\frac1{e(P_r)},
\qquad
\sigma\in\operatorname{Lin}(P_r).
\qquad\text{(9)}
\]

It is important that this is a measure on sequences. Permutations of
independent events are not identified.

\subsection{2.4. Restriction map}\label{restriction-map}

Since \(E_3\subset E_4\), define

\[
\rho_{43}:\operatorname{Lin}(P_4)\longrightarrow\operatorname{Lin}(P_3)
\qquad\text{(10)}
\]

by deleting all events in \(E_4\setminus E_3\) while preserving the
relative order of the remaining events. The suborder on \(E_3\) induced
by \(P_4\) coincides with \(P_3\), so the map is well defined. For every
new event \(N(a,b)\in E_4\setminus E_3\), both children belong to
\(T_3\), and hence all its prerequisites lie in \(E_3\). At the same
time, no new event can be a prerequisite of another new event, because
the children of level-\(4\) events belong to \(T_3\). Therefore
\(E_4\setminus E_3\) is an antichain of maximal elements, and every
predecessor of a new event lies in \(E_3\). Hence any sequence
\(\sigma\in\operatorname{Lin}(P_3)\) can be taken as the beginning of a
linear extension of \(P_4\), after which all new events can be appended
in any order. Thus \(\rho_{43}\) is surjective.

For \(\sigma\in\operatorname{Lin}(P_3)\), define

\[
F(\sigma)=|\rho_{43}^{-1}(\sigma)|.
\qquad\text{(11)}
\]

Then

\[
((\rho_{43})_*\mu_4)(\{\sigma\})
=\frac{F(\sigma)}{e(P_4)}.
\qquad\text{(12)}
\]

The pushforward is therefore uniform if and only if \(F(\sigma)\) is the
same for all lower-level orders.

The preceding restrictions are consistent. The map from level \(2\) to
level \(1\) is trivial because \(P_1\) has one event. For level \(3\) to
level \(2\), every linear extension of \(P_2\) begins with the unique
event \(A=N(L,L)\); the remaining three level-\(2\) events form an
antichain. After the first event, the number of released new events is
therefore the same for every lower order. Equivalently, the
release count is \((i+1)^2-4\) for every prefix length \(i\ge1\),
independent of the ordering of those three events. The recurrence in
Section 6 then gives equal fibers. Hence \(4\to3\) is the first
restriction in this grammar at which uniformity fails.

\section{3. Two witness orders}\label{two-witness-orders}

Set

\[
A=N(L,L),\quad
B=N(L,A),\quad
C=N(A,L),\quad
D=N(A,A).
\qquad\text{(13)}
\]

Then \(T_2=\{L,A,B,C,D\}\), and \(E_3\) consists of the \(25\) distinct
terms \(N(x,y)\) with \(x,y\in T_2\). To choose uniquely among several
available events, introduce the canonical string

\[
s(L)=\texttt{L},
\qquad
s(N(a,b))=\texttt{N(}s(a)\texttt{,}s(b)\texttt{)}.
\qquad\text{(14)}
\]

Define two linear extensions:

\begingroup
\raggedright
\begin{enumerate}
\def\labelenumi{\arabic{enumi}.}
\tightlist
\item
  \textbf{Level order:} at each step, choose the available event with
  the smallest key \((h(t),s(t))\).
\item
  \textbf{Depth-priority order:} at each step, choose the available
  event with the smallest key \((-h(t),s(t))\).
\end{enumerate}
\endgroup

Both rules select only an event whose prerequisites have already been
listed, and therefore generate linear extensions of \(P_3\). Denote them
by \(\sigma^{\mathrm{lev}}\) and \(\sigma^{\mathrm{dep}}\).

The positions of the four terms in \(E_2\) are decisive for the count
below:

\[
\begin{aligned}
\operatorname{pos}_{\mathrm{lev}}(A,B,C,D)&=(1,2,3,4),\\
\operatorname{pos}_{\mathrm{dep}}(A,B,C,D)&=(1,2,8,16).
\end{aligned}
\qquad\text{(15)}
\]

The first equality follows immediately from minimizing the height. For
the second, observe that after \(A\) the depth-priority rule selects
\(B\); after \(B\), it first exhausts five available terms of height
\(3\) and then selects \(C\); after \(C\), it exhausts seven more such
terms and only then selects \(D\). The full sequences are included in
\texttt{PROJECTIVE\_WITNESS.json} in the supplementary package.

\section{4. Main theorem}\label{main-theorem}

\textbf{Theorem 1.} For the grammar (1), partial orders (7), uniform
measures (9), and map (10),

\[
(\rho_{43})_*\mu_4\ne\mu_3.
\qquad\text{(16)}
\]

In particular, the two orders from Section 3 have fibers of unequal
cardinality:

\[
F(\sigma^{\mathrm{dep}})>F(\sigma^{\mathrm{lev}}),
\qquad\text{(17)}
\]

and the exact reduced ratio is

\[
\boxed{
\frac{F(\sigma^{\mathrm{dep}})}{F(\sigma^{\mathrm{lev}})}
=
\frac{614690215260160000}{479048686862260621}
}.
\qquad\text{(18)}
\]

We first prove the sign of the inequality without large computations. We
then describe two exact methods for obtaining (18).

\section{5. Analytic proof}\label{analytic-proof}

\subsection{5.1. The new events form an
antichain}\label{the-new-events-form-an-antichain}

From (4)--(5),

\[
|E_3|=25,
\qquad
|E_4\setminus E_3|=651.
\qquad\text{(19)}
\]

Let \(U=E_4\setminus E_3\). Every event \(u\in U\) has the form
\(N(a,b)\), where \(a,b\in T_3\), and at least one child does not belong
to \(T_2\). All its immediate prerequisites belong to \(E_3\). No event
in \(U\) is a prerequisite of another event in \(U\), because every
level-\(4\) event has children in \(T_3\). Thus \(U\) is an antichain
of maximal elements, and every predecessor of an element of \(U\) lies
in \(E_3\).

\subsection{5.2. Release profile}\label{release-profile}

Fix \(\sigma=(\sigma_1,\ldots,\sigma_{25})\in\operatorname{Lin}(P_3)\).
For \(u=N(a,b)\in U\), define the release index

\[
\tau_\sigma(u)
=
\max\bigl(
\{\operatorname{pos}_\sigma(a):a\ne L\}
\cup
\{\operatorname{pos}_\sigma(b):b\ne L\}
\bigr).
\qquad\text{(20)}
\]

The event \(u\) can be inserted into an extension immediately after the
\(\tau_\sigma(u)\)-th lower event has appeared. Let

\[
r_i^\sigma=|\{u\in U:\tau_\sigma(u)=i\}|,
\qquad
C_i^\sigma=\sum_{j=0}^{i}r_j^\sigma.
\qquad\text{(21)}
\]

Here \(r_0^\sigma=C_0^\sigma=0\): no event in \(U\) can have both
children equal to the leaf, because that term already belongs to
\(E_1\subset E_3\).

The quantity \(C_i^\sigma\) is the cumulative release profile: the
number of new events available after the first \(i\) lower events.

For a prefix, define

\[
S_i^\sigma=\{L,\sigma_1,\ldots,\sigma_i\},
\qquad
b_i^\sigma=|S_i^\sigma\cap T_2|.
\qquad\text{(22)}
\]

The set \(S_i^\sigma\) has \(i+1\) elements. Ordered pairs from it
specify \((i+1)^2\) available terms of the form \(N(a,b)\). Exactly
\((b_i^\sigma)^2\) pairs have both children in \(T_2\) and therefore
already belong to \(E_3\) rather than \(U\). Hence

\[
\boxed{
C_i^\sigma=(i+1)^2-(b_i^\sigma)^2
}.
\qquad\text{(23)}
\]

This is an exact formula, not an asymptotic approximation. At \(i=25\),
it gives \(26^2-5^2=651\).

\subsection{5.3. Comparison of the two
profiles}\label{comparison-of-the-two-profiles}

The positions in (15) give

\[
b_i^{\mathrm{lev}}=
\begin{cases}
1+i,&0\le i\le4,\\
5,&5\le i\le25,
\end{cases}
\qquad\text{(24)}
\]

and

\[
b_i^{\mathrm{dep}}=
\begin{cases}
1,&i=0,\\
2,&i=1,\\
3,&2\le i\le7,\\
4,&8\le i\le15,\\
5,&16\le i\le25.
\end{cases}
\qquad\text{(25)}
\]

Substitution into (23) yields

\[
C_i^{\mathrm{dep}}-C_i^{\mathrm{lev}}=
\begin{cases}
0,&0\le i\le2,\\
7,&i=3,\\
16,&4\le i\le7,\\
9,&8\le i\le15,\\
0,&16\le i\le25.
\end{cases}
\qquad\text{(26)}
\]

Thus the depth-priority profile is nowhere smaller and is strictly
larger on a nonempty interval of prefix lengths.

\subsection{5.4. Monotonicity of the number of
interleavings}\label{monotonicity-of-the-number-of-interleavings}

\textbf{Lemma 1.} Let a fixed lower chain contain \(n\) elements, with
\(m\) pairwise incomparable labeled upper events above it. If two
cumulative profiles satisfy \(C'_i\ge C_i\) for every \(0\le i\le n\),
then the number of admissible interleavings for \(C'\) is at least the
number for \(C\). If the inequality is strict for at least one \(i<n\),
then the number of interleavings is strictly larger.

\textbf{Proof.} Write the release indices of the upper events in
nondecreasing order:

\[
\tau_1\le\cdots\le\tau_m,
\qquad
\tau'_1\le\cdots\le\tau'_m.
\]

Pointwise dominance of the cumulative profiles is equivalent to
\(\tau'_k\le\tau_k\) for every \(k\). Label the corresponding upper
events \(u_k\) and \(u'_k\), thereby defining a bijection
\(u_k\mapsto u'_k\).

Take an admissible interleaving for the first profile and replace each
label \(u_k\) by \(u'_k\) without changing its position. The event
\(u_k\) occurred after the \(\tau_k\)-th lower element, and
\(\tau'_k\le\tau_k\); hence the new word is admissible. The relabeling
is invertible on its image, so this defines an injection.

If \(\tau'_k<\tau_k\) for some \(k\), place \(u'_k\) immediately after
the \(\tau'_k\)-th lower element, place the rest of the lower chain
after it, and put all other upper events at the end. If \(\tau'_k=0\),
this means placing \(u'_k\) at the beginning of the word, before the
first lower element. This word is admissible for the second profile but
is not in the image of the injection: its putative preimage would
contain \(u_k\) before its release. Thus the injection is not
surjective. \(\square\)

After the upper events have been matched in nondecreasing order of
release index, the poset associated with the dominating profile is
obtained from the other one by deleting comparabilities between the
lower chain and matched upper events. Thus the weak part of Lemma 1 is
also a special case of the standard monotonicity of the number of linear
extensions under deletion of comparabilities {[}14, Theorem 3.8{]}. The
proof above is retained because it gives the injection explicitly and
proves strictness directly for the two fibers used here.

\subsection{5.5. Completion of the inequality
proof}\label{completion-of-the-inequality-proof}

Once \(\sigma\) is fixed, the internal dependencies of \(P_3\) are
already encoded by the lower chain, and every extension to \(P_4\) is an
admissible interleaving of this chain with \(651\) labeled upper events.
Equation (26) and Lemma 1 give

\[
F(\sigma^{\mathrm{dep}})>F(\sigma^{\mathrm{lev}}).
\qquad\text{(27)}
\]

By (12), two points in \(\operatorname{Lin}(P_3)\) have different
probabilities under \((\rho_{43})_*\mu_4\), whereas they are
equiprobable under \(\mu_3\). This proves (16) without using large
integers.

\section{6. Exact sizes of the two
fibers}\label{exact-sizes-of-the-two-fibers}

The ratio (18) was obtained from two different decompositions of the
same finite set of extensions. Both use integer arithmetic.

\subsection{6.1. Counting by available
events}\label{counting-by-available-events}

Let \(n=25\), \(m=651\), and let \(C_i\) be the cumulative profile of a
fixed lower order. Denote by \(Q(i,j)\) the number of admissible
completions of a state in which the first \(i\) lower events and \(j\)
upper events have already been listed. Then

\[
Q(n,m)=1,
\qquad\text{(28)}
\]

and for all other states,

\[
Q(i,j)=
\mathbf 1_{i<n}Q(i+1,j)
+
\mathbf 1_{j<C_i}(C_i-j)Q(i,j+1).
\qquad\text{(29)}
\]

The first term appends the uniquely determined next lower event. The
second chooses one of the \(C_i-j\) upper labels that have been released
but not yet used. Every complete word has a unique sequence of such
states, and therefore

\[
F(\sigma)=Q(0,0).
\qquad\text{(30)}
\]

\subsection{6.2. Counting by positions of the lower
chain}\label{counting-by-positions-of-the-lower-chain}

The second method fixes the positions of the lower events in the
complete word of length \(n+m\):

\[
1\le p_1<\cdots<p_n\le n+m.
\qquad\text{(31)}
\]

Let

\[
r_i=C_i-C_{i-1},
\qquad
H_i=\sum_{k=i+1}^{n}r_k,
\qquad\text{(32)}
\]

and define the falling factorial by

\[
(q)_{\underline{k}}=
\begin{cases}
q(q-1)\cdots(q-k+1),&q\ge k\ge0,\\
0,&q<k\ \text{or}\ k<0.
\end{cases}
\qquad\text{(33)}
\]

After position \(p_i\), there are \(m+i-p_i\) positions not occupied by
the lower chain. If upper events are assigned in descending order of
threshold, then \(H_i\) positions have already been occupied by events
with a later threshold. The \(r_i\) distinct events with exact threshold
\(i\) can therefore be placed in

\[
(m+i-p_i-H_i)_{\underline{r_i}}
\]

ways. Summing over all positions of the lower chain gives

\[
\boxed{
F(\sigma)=
\sum_{1\le p_1<\cdots<p_n\le n+m}
\prod_{i=1}^{n}
(m+i-p_i-H_i)_{\underline{r_i}}
}.
\qquad\text{(34)}
\]

Equation (34) counts every admissible word exactly once: the word
uniquely determines the positions of the lower chain, and the remaining
labeled positions are then assigned uniquely among the threshold groups.

The sum need not be evaluated by enumerating all \(\binom{n+m}{n}\)
position sets. If

\[
G_1(p)=(m+1-p-H_1)_{\underline{r_1}},
\]

then the recurrence

\[
G_i(p)=
(m+i-p-H_i)_{\underline{r_i}}
\sum_{q<p}G_{i-1}(q)
\qquad\text{(35)}
\]

gives

\[
F(\sigma)=\sum_{p=1}^{n+m}G_n(p).
\qquad\text{(36)}
\]

The first method classifies prefixes by the numbers of lower and upper
events already listed. The second classifies complete words by the
global positions of the lower chain. They use the same grammar and
profiles, but different state spaces and recurrences.

\subsection{6.3. Counting result}\label{counting-result}

Both methods produced the same value for each fiber.

\begin{longtable}[]{@{}
  >{\raggedright\arraybackslash}p{(\linewidth - 4\tabcolsep) * \real{0.2857}}
  >{\centering\arraybackslash}p{(\linewidth - 4\tabcolsep) * \real{0.3571}}
  >{\centering\arraybackslash}p{(\linewidth - 4\tabcolsep) * \real{0.3571}}@{}}
\caption{Exact check values for the two fibers.}\tabularnewline
\toprule\noalign{}
\begin{minipage}[b]{\linewidth}\raggedright
Quantity
\end{minipage} & \begin{minipage}[b]{\linewidth}\centering
Level order
\end{minipage} & \begin{minipage}[b]{\linewidth}\centering
Depth-priority order
\end{minipage} \\
\midrule\noalign{}
\endfirsthead
\toprule\noalign{}
\begin{minipage}[b]{\linewidth}\raggedright
Quantity
\end{minipage} & \begin{minipage}[b]{\linewidth}\centering
Level order
\end{minipage} & \begin{minipage}[b]{\linewidth}\centering
Depth-priority order
\end{minipage} \\
\midrule\noalign{}
\endhead
\bottomrule\noalign{}
\endlastfoot
Events in \(E_3\) & 25 & 25 \\
New events in \(E_4\setminus E_3\) & 651 & 651 \\
Sum of \(r_i\) & 651 & 651 \\
Decimal digits in \(F(\sigma)\) & 1557 & 1557 \\
\(\log_{10}F\), for orientation only & 1556.259117544793456 &
1556.367394190755243 \\
SHA-256, abbreviated & \texttt{73ABE791...94CA04E} &
\texttt{56D0E809...86BB88E8} \\
\end{longtable}

The main JSON certificate contains the complete decimal expansions and
checksums, and the CSV table duplicates both integers. Reducing their
ratio by the greatest common divisor gives (18); the numerator and
denominator are coprime. Neither the decimal logarithms nor the
approximation enter the verification.

\section{7. Correctness checks and
reproduction}\label{correctness-checks-and-reproduction}

The short analytic proof of the inequality does not depend on program
code. The computational part is needed for the exact value (18) and is
accompanied by the following checks:

\begin{itemize}
\item
  explicit generation reproduces \(|T_0|,\ldots,|T_4|=1,2,5,26,677\);
\item
  both orders contain exactly all \(25\) events in \(E_3\) and respect
  all prerequisites;
\item
  the profiles obtained directly from the children of the terms agree
  with (23);
\item
  each histogram \((r_0,\ldots,r_{25})\) sums to \(651\);
\item
  methods (29) and (34)--(36) agree separately for both fibers;
\item
  an independent implementation of the second method does not import the
  main verifier and reproduces both integers and the fraction;
\item
  the complete SHA-256 hashes of the decimal expansions are recorded in
  the certificate, and the prefixes and suffixes shown in the table
  agree with them;
\item
  every equality test is carried out over integers and rational numbers.
\end{itemize}

The supplementary directory contains:

\begin{Shaded}
\begin{Highlighting}[]
\NormalTok{anc/}
\NormalTok{\textasciigrave{}{-}{-} 04\_LINEAR\_EXTENSION\_NONPROJECTIVITY/}
\NormalTok{    |{-}{-} ARTICLE\_SOURCE\_EN.md}
\NormalTok{    |{-}{-} CITATION.txt}
\NormalTok{    |{-}{-} verify\_linear\_extension\_nonprojectivity.py}
\NormalTok{    |{-}{-} verify\_independent\_projective\_witness.py}
\NormalTok{    |{-}{-} README\_EN.md}
\NormalTok{    |{-}{-} README\_UK.md}
\NormalTok{    |{-}{-} requirements.txt}
\NormalTok{    |{-}{-} SHA256SUMS.txt}
\NormalTok{    \textasciigrave{}{-}{-} RUN/}
\NormalTok{        |{-}{-} LINEAR\_EXTENSION\_NONPROJECTIVITY\_CERTIFICATE.json}
\NormalTok{        |{-}{-} PROJECTIVE\_WITNESS.json}
\NormalTok{        |{-}{-} ORDER\_EXTENSION\_COUNTS.csv}
\NormalTok{        \textasciigrave{}{-}{-} INDEPENDENT\_CERTIFICATE.json}
\end{Highlighting}
\end{Shaded}

Both scripts use only the Python standard library. The main script
performs \(25\) exact checks; the field \texttt{all\_checks\_pass} in
its certificate must have the value \texttt{true}. The independent
script must terminate with the field \texttt{pass:\ true}. The files use
relative paths and create the results directory themselves.

As an additional small check, we computed

\[
e(P_1)=1,
\qquad
e(P_2)=6,
\qquad
e(P_3)=861733891296165888000.
\qquad\text{(37)}
\]

All six orders of \(P_2\) have the same number of extensions to \(P_3\).
Thus, among the four directly tested restrictions from levels \(1\),
\(2\), \(3\), and \(4\) to the preceding level, the earliest failure in
this construction occurs in the transition from level \(4\) to level
\(3\). This statement does not extend to other families or other maps.

\section{8. Relation to the
literature}\label{relation-to-the-literature}

A uniformly random linear extension of a finite partially ordered set is
a classical object; for a recent survey, see {[}14{]}. Exact counting of linear extensions in general is a
\(\#P\)-complete problem {[}5{]}, while near-uniform sampling can be
constructed using specially designed rapidly mixing Markov chains
{[}6{]}. Our problem is not a general counting algorithm: the special
two-level structure of the fiber reduces it to a threshold profile.

Order-invariant measures on natural extensions of causal sets concern
conditional equiprobability of admissible orders and existence and
uniqueness questions for probability measures on infinite natural
extensions {[}1,2{]}. Central measures on paths of
graded graphs are described through coherent cotransition data and
projective limits of simplices {[}3,4,11{]}. For total vertex orderings
of hereditary graph classes, compatibility under passage to induced
substructures has also been studied {[}7,8{]}. Sampling consistency for
rooted binary trees provides a further nearby analogy, but its
projection restricts a tree to a selected subset of leaves and
suppresses degree-two vertices other than the root, rather than restricting a
linear extension of a prerequisite poset {[}13{]}. These works establish the
broader setting: uniformity on an individual slice and coherence of a
family are different requirements.

The counterexample given here is not a new general theory of central
measures, nor is it the first example of nonuniform fibers under
deletion. Its specificity lies in the exact conjunction of four
ingredients: the nested grammar (1), consistency through the preceding
levels, natural deletion of \(E_4\setminus E_3\), and the uniform
measure on linear extensions. To the best of our knowledge, this delayed
level-\(4\) to level-\(3\) failure, its release-profile injection, and
the exact ratio (18) have not been treated previously. This is a limited
novelty statement, not a claim of absolute priority.

Models of classical sequential growth of causal sets impose causality
conditions and a discrete form of general covariance on transition
probabilities {[}12{]}. Our result
is not such a growth model and has no independent physical
interpretation. It shows only, through a specific finite example, why
separate uniform normalization at every level cannot replace a check of
interlevel consistency.

\section{9. Scope of the claim and open
problems}\label{scope-of-the-claim-and-open-problems}

Theorem 1 depends on all of the following conditions.

\begin{enumerate}
\def\labelenumi{\arabic{enumi}.}
\tightlist
\item
  \textbf{Exact identity of events.} Events are syntactic terms.
  Quotienting by tree isomorphisms or introducing separately labeled
  copies changes the state space.
\item
  \textbf{Ordered constructor.} The pairs \((a,b)\) and \((b,a)\) are
  distinct. An unordered constructor changes the recurrence and the
  profiles.
\item
  \textbf{Conflict-free persistence.} Events are not consumed, and the
  new layer is an antichain. Conflicts or dependencies within the new
  layer invalidate the reduction to (23).
\item
  \textbf{Uniform weights.} A different coherent family of measures may
  exist. The counterexample does not exclude nonuniform cotransition
  weights.
\item
  \textbf{Specified restriction.} The map used here deletes the new
  labels. A different deterministic or stochastic coarse graining has
  different fibers.
\item
  \textbf{Sequences as states.} If all permutations of independent
  events are identified, only one complete unordered configuration
  remains at each level, and the randomness studied here disappears.
\end{enumerate}

The theorem uses only two fibers. This is sufficient to disprove
uniformity, but not to describe the entire distribution

\[
\{F(\sigma):\sigma\in\operatorname{Lin}(P_3)\}.
\]

The minimum and maximum fiber sizes, the number of distinct values, the
variance of the pushforward, the asymptotic behavior as \(r\) grows, and
the classification of weights that form a projective system on this
grammar remain open. A separate natural question is whether such weights
can be described by a local rule that does not require a new set of
parameters at every level.

\section{10. Conclusions}\label{conclusions}

For the nested grammar of ordered binary terms, the uniform measure on
linear extensions is well defined at each individual finite level, but
the resulting family of measures is not projectively consistent. The
restrictions through level \(3\) are consistent; at the first failure,
from level \(4\) to level \(3\), two explicit lower-level orders have
different numbers of extensions.

The source of the inequality is simple once the problem is expressed in
the right variables. Each lower prefix releases \((i+1)^2-b_i^2\) new
events. The depth-priority order introduces more deeply nested terms
earlier and therefore has strictly more available upper insertions over
an interval of prefix lengths. The injection lemma turns this local
advantage into a strict global result.

The exact magnitude of the effect is confirmed by two different integer
recurrences. The ratio of the fiber cardinalities is (18), and the
complete integers and machine certificates are included in the
supplementary package. The strength of the conclusion lies not in
universality but in complete specification: the grammar, map, two
witnesses, analytic sign, and exact magnitude are all given without
hidden parameters.

\section{Declarations}\label{declarations}

\textbf{Funding.} This work received no external funding.

\textbf{Conflict of interest.} The author declares no conflict of
interest.

\textbf{Author contribution.} The sole author formulated the problem,
proved the structural results, implemented and checked the exact
enumerations, and approved the manuscript.

\textbf{AI-assisted tools.} Generative-AI tools were used during
drafting, language editing, code review, and computational
cross-checking. The author independently checked the proofs, references,
and executable artifacts and assumes full responsibility for the
contents.

\textbf{Ethics approval.} Not applicable: this work includes no research
involving human participants or animals.

\textbf{Data and code availability.} All data are finite synthetic
objects generated by definition (1). The code and exact certificates are
included in \texttt{anc/04\_LINEAR\_EXTENSION\_NONPROJECTIVITY/} in the
arXiv source package.

\section{References}\label{references}

\begin{enumerate}
\def\labelenumi{\arabic{enumi}.}
\tightlist
\item
  G. Brightwell, M. Łuczak. ``Order-invariant measures on causal sets.''
  \emph{The Annals of Applied Probability} \textbf{21}(4), 1493-1536
  (2011). \url{https://doi.org/10.1214/10-AAP736}
\item
  G. Brightwell, M. Łuczak. ``Order-invariant measures on fixed causal
  sets.'' \emph{Combinatorics, Probability and Computing}
  \textbf{21}(3), 330-357 (2012).
  \url{https://doi.org/10.1017/S0963548311000721}
\item
  A. M. Vershik. ``The Problem of Describing Central Measures on the
  Path Spaces of Graded Graphs.'' \emph{Functional Analysis and Its
  Applications} \textbf{48}(4), 256-271 (2014).
  \url{https://doi.org/10.1007/s10688-014-0069-5}
\item
  A. M. Vershik. ``Equipped Graded Graphs, Projective Limits of
  Simplices, and Their Boundaries.'' \emph{Journal of Mathematical
  Sciences} \textbf{209}(6), 860-873 (2015).
  \url{https://doi.org/10.1007/s10958-015-2533-z}
\item
  G. Brightwell, P. Winkler. ``Counting linear extensions.''
  \emph{Order} \textbf{8}(3), 225-242 (1991).
  \url{https://doi.org/10.1007/BF00383444}
\item
  R. Bubley, M. Dyer. ``Faster Random Generation of Linear Extensions.''
  \emph{Discrete Mathematics} \textbf{201}(1-3), 81-88 (1999).
  \url{https://doi.org/10.1016/S0012-365X(98)00333-1}
\item
  O. Angel, A. S. Kechris, R. Lyons. ``Random Orderings and Unique
  Ergodicity of Automorphism Groups.'' \emph{Journal of the European
  Mathematical Society} \textbf{16}(10), 2059-2095 (2014).
  \url{https://doi.org/10.4171/JEMS/483}
\item
  P. N. Balister, B. Bollobás, S. Janson. ``Consistent Random
  Vertex-Orderings of Graphs.'' \emph{Journal of the European
  Mathematical Society} \textbf{27}(7), 2623-2652 (2025).
  \url{https://doi.org/10.4171/JEMS/1624}
\item
  A. V. Aho, N. J. A. Sloane. ``Some Doubly Exponential Sequences.''
  \emph{The Fibonacci Quarterly} \textbf{11}(4), 429-437 (1973).
  \url{https://doi.org/10.1080/00150517.1973.12430815}
\item
  P. Flajolet, A. M. Odlyzko. ``Limit Distributions for Coefficients of
  Iterates of Polynomials with Applications to Combinatorial
  Enumerations.'' \emph{Mathematical Proceedings of the Cambridge
  Philosophical Society} \textbf{96}(2), 237-253 (1984).
  \url{https://doi.org/10.1017/S0305004100062149}
\item
  A. M. Vershik. ``One-dimensional central measures on numberings of
  ordered sets.'' \emph{Functional Analysis and Its Applications}
  \textbf{56}(4), 251-256 (2022).
  \url{https://doi.org/10.1134/S0016266322040025}
\item
  D. P. Rideout, R. D. Sorkin. ``Classical Sequential Growth Dynamics
  for Causal Sets.'' \emph{Physical Review D} \textbf{61}, 024002
  (2000). \url{https://doi.org/10.1103/PhysRevD.61.024002}
\item
  B. Hollering, S. Sullivant. ``Exchangeable and Sampling-Consistent
  Distributions on Rooted Binary Trees.'' \emph{Journal of Applied
  Probability} \textbf{59}(1), 60--80 (2022).
  \url{https://doi.org/10.1017/jpr.2021.28}
\item
  S. H. Chan, I. Pak. ``Linear Extensions of Finite Posets.''
  \emph{EMS Surveys in Mathematical Sciences} (2025), published online
  first. \url{https://doi.org/10.4171/EMSS/97}
\end{enumerate}

\end{document}